\documentclass[12pt,a4paper]{article}
\usepackage{amsmath}
\usepackage{amssymb}
\usepackage{amsthm}
\usepackage{latexsym}
\usepackage[english,russian]{babel}

\newcommand{\Ker}{\mathop{\rm Ker}\nolimits}
\newcommand{\im}{\mathop{\rm Im}\nolimits}

\newcommand{\tr}{\mathop{\rm tr}\nolimits}
\newcommand{\sgn}{\mathop{\rm sgn}\nolimits}
\newcommand{\vol}{\mathop{\rm Vol}\nolimits}

\newcommand{\Lvar}{\mathop{\rm Loc.var}\nolimits}

\newcommand{\abs}[1]{\lvert#1\rvert}

\numberwithin{equation}{section}
\swapnumbers

\newtheorem{Theorem}{Теорема}[section]
\newtheorem{Proposition}{Предложение}[section]

\newtheorem{Lemma}[Theorem]{Лемма}

\theoremstyle{definition}
\newtheorem{defn}{Определение}[section]
\theoremstyle{remark}

\def\?{\marginpar{$\bullet\bullet\bullet$}}

\title{О параллелепипедах в алгебре и топологии}
\author{Баяк И. В.}
\date{23 марта 2007}

\begin{document}
\maketitle
\begin{abstract}

Концепцию применения параллелепипедов можно найти в таких разделах
математики как комбинаторный анализ, полилинейная алгебра и
алгебраическая топология. Действительно, последовательности ребер
n--мерного параллелепипеда, натянутого на базис n--мерного
пространства, изоморфны подстановкам. В свою очередь, вычисление
ориентированного объема параллелепипеда как полилинейной функции,
принимающей нулевое значение на линейно зависимых векторах,
приводит к понятию знакопеременного полилинейного произведения.
Наконец, параллелепипеды легко могут складываться в клеточные
пространства и образовывать цепные комплексы с группой гомологий в
качестве топологического инварианта этих пространств. Имея в виду
вышеперечисленные проникновения геометрии в алгебру и топологию,
мы будем здесь последовательно развивать алгебро-топологический
формализм, имеющий отношение к параллелепипедам, а затем найдем
ему применение в решении проблемы классификации замкнутых
ориентируемых многообразий произвольной размерности.

\end{abstract}

\section{Подстановки}

Прежде всего обратимся к основаниям комбинаторного анализа, т. е.
к понятию отображения множеств, где устоялись классические
определения, и поэтому нам будет достаточно сослаться на
университетский учебник, например, \cite{тыш}. Однако, даже
классические представления иногда допускают определенную
модернизацию. Действительно, пусть дана такая тройка множеств
$(A,B,F)$, что $F=\{(a,b)\}\ne \varnothing$, где $a\in A$, $b\in
B$. Тогда мы определим неклассическое понятие отображения $F$ с
классическим обозначением $f:A \to B$, имея в виду, что каждая
пара из $F$ может быть представлена стрелкой соответствия $a \to
b$, где $b=f(a)$ для однозначного соответствия и $b \in f(a)$ для
многозначного соответствия. Вместе с тем, триплет $(B,A,F^{-1})$,
где имеет место обращение всех стрелок, задает обратное
отображение. Если $\forall \: a\: \exists \: b \quad (f(a)=b)$, то
мы говорим об однозначном отображении всего $A$ и тем самым
возвращаемся к классическому представлению об отображении.
Сужением отображения $F$ на $A'\subset A$ мы назовем такое
отображение $f'=f\vert_{A'} :A'\to B$, которое соответствует
триплету $(A',B,F')$, где включение $F'\subset F$ порождается
включением $A'\subset A$; т.е., исключением из $F$ стрелок,
выходящих из подмножества $A\backslash A'$.

Расширение понятия отображения влечет за собой расширение таких
понятий как инъективное и сюръективное отображение. Инъективным
отображением можно было бы назвать такое $F$, которое своим первым
аргументом исчерпывает $A$, и в котором нет сходящихся стрелок, т.
е. отсутствуют пары $(x,b),(y,b)$, иначе говоря, $\forall \: x\ne
y\in A\Rightarrow f(x)\cap f(y)=\varnothing $. В свою очередь,
сюръективным отображением можно было бы назвать такое $F$, которое
своим вторым аргументом исчерпывает $B$, и в котором нет
расходящихся стрелок, т. е. отсутствуют пары $(a,x),(a,y)$, иначе
говоря, $\forall \: b \: \exists \: a \quad (f(a)=b)$.
Квазиинъективным отображением можно назвать такое $F$, которое
исчерпывает $A$, и в котором нет расходящихся стрелок, иначе
говоря, $\forall \: a \: \exists \: b \quad (f(a)=b)$.
Квазисюръективным отображением можно назвать такое $F$, которое
исчерпывает $B$, и в котором нет сходящихся стрелок, иначе говоря,
$\forall \: b \: \exists ! \: a \quad (f(a)=b)$. Наконец, строго
инъективным отображением можно назвать отображение, инъективное в
классическом смысле, т. е. такое, что $\forall \: x\ne y\in
A\Rightarrow f(x)\neq f(y)$. Аналогично, строго сюръективным
отображением можно назвать отображение, сюръективное в
классическом смысле. Вместе с тем, в понятие биективного
отображения мы, по--прежнему, вкладываем классический смысл.

Логические связи между всеми этими типами отображений отражены в
нижеследующих очевидных утверждениях.

\begin{Proposition}

Биективное отображение одновременно инъективно и сюръективно

\end{Proposition}

\begin{Proposition}

Сужение биекции строго инъективно

\end{Proposition}

\begin{Proposition}

Отображение, обратное инъективному (квазиинъективному), является
сюръекцией (квазисюръекцией)

\end{Proposition}

Теперь, в рамках вышепринятого формализма, нам будет легко
определить основные комбинаторные понятия. Действительно, если
произвольную биекцию из триплета $(A,A,F)$ называть подстановкой,
а простейшую подстановку --- транспозицией, то размещение без
повторений можно определить как строгую инъекцию, а размещение с
повторениями --- как квазиинъекцию. Примем также удобную нам
систему обозначений. Если заданы множества  $I= \{1, \ldots ,n \},
\quad J= \{1, \ldots ,m \}$, где $m \leq n$, тогда отображение
$/J/:=/i_{1}, \ldots ,i_{j}, \ldots ,i_{m}/:J \to J$ пусть
обозначает произвольную подстановку $m$ элементов множества $J$, а
$\sgn /J/$ --- знак четности этой подстановки; отображение
$[J]:=[i_{1}, \ldots ,i_{j}, \ldots ,i_{m}]:J \rightarrow I$ пусть
обозначает произвольное размещение из $n$ по $m$ с повторением.
Отображение $\langle J \rangle :=\langle i_{1}, \ldots ,i_{j},
\ldots ,i_{m} \rangle:J \rightarrow I$ обозначает произвольное
размещение из $n$ по $m$ без повторений; отображение $(J):=(i_{1},
\ldots ,i_{j}, \ldots ,i_{m}):J \rightarrow I$ пусть обозначает
произвольное упорядоченное размещение, т. е. такое $\langle J
\rangle$ в котором $i_{j}<i_{j+1}$. Наконец, если множество $\{ J
\} :=\{ i_{1}, \ldots ,i_{j}, \ldots ,i_{m} \}$ обозначает у нас
производное от $(J)$ подмножество $I$, состоящее из $m$ элементов
образа размещения $(J)$, тогда отображение $/(J)/:=/(i_{1}, \ldots
,i_{j}, \ldots ,i_{m})/:\{ J \}\rightarrow \{ J \}$ пусть
обозначает произвольную подстановку элементов множества $\{ J \}$,
соответственно $\sgn /(J)/$ --- знак четности этой подстановки.

Пусть отображение $/(J)\hat{i}_{j}/: (i_{1},\ldots,
i_{j},\ldots,i_{m})\rightarrow (i_{1},\ldots,
\hat{i}_{j},\ldots,i_{m})i_{j}$ обозначает производную от $(J)$ и
произвольного элемента $\{J\}$ подстановку изъятия этого элемента,
в результате которой элемент $i_{j}$ переходит в конец а остальные
элементы, после переноса $i_{j}$, остаются упорядоченными по
возрастанию. Наряду с изъятием элемента из упорядоченного
множества можно рассматривать присоединение элемента с последующим
упорядочиванием. Тогда мы получим подстановку $/(J)\check{k}/$,
где $k\in \{J'\}= I\setminus \{J\}$, которую следует понимать как
подстановку упорядочивания множества $\{(J),k\}$. Все эти
подстановки обладают свойством антикоммутативности их четностей,
так что $\sgn /(J)\hat{i}_{j}\hat{i}_{k}/= -\sgn
/(J)\hat{i}_{k}\hat{i}_{j}/$, а $\sgn /(J)\check{j}\check{k}/=
-\sgn /(J)\check{k}\check{j}/$.

Рассмотрим, наконец, подстановки частичного нарушения
упорядоченности множества $I$. Принимаем, что образ подстановки
$/J,J'/: I\rightarrow \{(J),(J')\}$, где $\{J'\}=I\setminus\{J\}$,
образован множествами $\{J\}$ и $\{J'\}$, упорядоченными по
отдельности, а подстановка $/(J,J')/: \{(J),(J')\}\rightarrow
\{(J'),(J)\}$ подразумевает перестановку упорядоченных множеств.
Тогда, поскольку $/J',J/= /J,J'/*/(J,J')/$ а $\sgn /(J,J')/=
(-1)^{m(n-m)}$, то мы получим равенство $\sgn /J,J'/\cdot\sgn
/J',J/=(-1)^{m(n-m)}$.

\section{Знакопеременные произведения}

Алгебра знакопеременных полилинейных произведений хорошо изучена,
и поэтому можно сослаться, например, на соответствующую главу в
\cite{лен} или \cite{ефи}. Здесь же мы, используя понятие четности
частично упорядоченных подстановок и подстановок изъятия
(присоединения), коснемся лишь отдельных вопросов.

Прежде всего обратимся к полилинейным произведениям в целом. Пусть
$\mathbb{L}$ означает $n$--мерное линейное пространство над
$\mathbb{R}$ а $(e^{i})_{I}$ --- его базис. Пусть также
$\mathbb{T}^{m}(\mathbb{L})$ означает $n^{m}$--мерное пространство
над $\mathbb{R}$ а $(e^{[J]})_{\{ [J] \}}$ --- его базис. Тогда мы
можем явно представить тензорное полилинейное произведение, если
запишим
$$x^{1}\otimes \cdots \otimes x^{m}:\mathbb{L}^{m}\rightarrow
\mathbb{T}^{m}(\mathbb{L}):(x^{1},\ldots ,x^{m})\rightarrow
\sum_{\{ [J] \}}\prod_{i_{j}\in [J]}x_{ij}e^{[J]},$$ где $x_{ij}$
означает $j$--ю координату вектора $x^{i}$. Поскольку
$e^{[i_{1}}\otimes\cdots\otimes e^{i_{m}]}=e^{[J]}$, то линейная
оболочка множества $\{ e^{[i_{1}}\otimes\cdots\otimes e^{i_{m}]}
\}_{\{ [J] \}}$ изоморфна тензорному линейному пространству
$\mathbb{T}^{m}(\mathbb{L})$.

Сформируем теперь два знакопеременных полилинейных произведения, а
именно, внешнее произведение $$x^{1}\wedge \cdots\wedge
x^{m}:\mathbb{L}^{m}\rightarrow \Lambda^{m}(\mathbb{L}):$$
$$(x^{1},\ldots ,x^{m})\rightarrow\sum_{\{ (J) \}}e^{(J)}\sum_{\{
/(J)/ \}}\sgn /(J)/\prod_{i_{j}\in /(J)/}x_{ij}=\sum_{\{ (J)
\}}\det(x_{ij})^{\{ J \}}_{J}e^{(J)},$$ где
$\Lambda^{m}(\mathbb{L})$ это пространство поливекторов,
изоморфное линейной оболочке множества
$\{e^{(i_{1}}\wedge\cdots\wedge e^{i_{m})}\}_{\{(J)\}}$; и
векторное произведение $$x^{1}\times\cdots\times x^{m}:
\mathbb{L}^{m} \rightarrow\Lambda^{n-m}(\mathbb{L}):$$
$$(x^{1},\ldots ,x^{m})\rightarrow\sum_{\{ (J')
\}}\sgn/J,J'/\det(x_{ij})^{\{ J \}}_{J}e^{(J')},$$ где $\{J'\}=
I\setminus\{J\}$. Всем этим произведениям гарантирована
знакопеременность, поскольку при транспозиции произвольной пары
векторов из набора $(x^{1},\ldots ,x^{m})$ знаки четностей всех
подстановок $/(J)/$ меняются на противоположные. Кроме того, для
них справедливо
\begin{Proposition}
Система векторов $(x^{1},\ldots ,x^{m})$ линейно зависима тогда и
только тогда, когда $x^{1}\wedge\cdots\wedge
x^{m}=x^{1}\times\cdots\times x^{m}=0$
\end{Proposition}
Действительно, с одной сторолны, если система $(x^{1},\ldots
,x^{m})$ линейно зависима, то $\det(x_{ij})^{\{ J
\}}_{J}=0\quad\forall\:(J)\in\{(J)\}$, с другой стороны, если
система $(x^{1},\ldots ,x^{m})$ линейно независима, то с помощью
элементарных преобразований системы и перестановок элементов
базиса $(e^{1},\ldots,e^{n})$ матрица $(x_{ij})^{I}_{J}$
приводится к трапецеидальному виду с $m$ ненулевыми строками и
ненулевой диагональю $(x_{ii})_{J}$, в силу чего
$\det(x_{ij})_{J}^{J}\neq 0$. Здесь мы воспользовались тем
обстоятельством, что элементарные преобразования сохраняют
независимость системы векторов и одновременно сохраняют равенство
(неравенство) нулю определителя, а подстановки базиса не зануляют
ни внешнего ни векторного произведения. Тем самым, внешнее и
векторное произведение линейно независимой системы векторов имеют
по крайней мере по одной ненулевой координате, а следовательно
утверждение доказано.

Пространство поливекторов $\Lambda^{m}(\mathbb{L})$ интересно еще
и тем, что для него существуют естественные понижающий и
повышающий граничные гомоморфизмы, а именно,
$$\delta:\Lambda^{m}(\mathbb{L})\rightarrow
\Lambda^{m-1}(\mathbb{L}):e^{(J)}\rightarrow
\sum_{j\in\{J\}}\sgn/(J)\hat{j}/\cdot e^{\hat{(J)}},$$ где
$\hat{(J)}$ означает упорядоченную подстановку изъятия
$/(J)\hat{j}/$, и $$d:\Lambda^{m}(\mathbb{L})\rightarrow
\Lambda^{m+1}(\mathbb{L}):e^{(J)}\rightarrow\sum_{j\in I\backslash
\{J\}}\sgn/(J)\check{j}/\cdot e^{\check{(J)}},$$ где $\check{(J)}$
означает упорядоченную подстановку присоединения $/(J)\check{j}/$.
Действительно, поскольку в поливекторе $\delta\delta(x)$ базисные
элементы образуются с помощью изъятий дважды, то в силу свойства
антикоммутативности четностей подстановок с изъятием мы получим
тождество $\delta\delta(x)=0$ для всякого поливектора $x$.
Аналогично получается тождество $dd(x)=0$.

Обратимся теперь к пространству $\mathbb{L}$, в котором задан
функционал скалярного произведения. Пусть $\mathbb{E}$ будет
произвольным n--мерным евклидовым пространством над $\mathbb{R}$.
Если принять, что $$x^{1}\wedge\cdots\wedge x^{m}\cdot
y^{1}\wedge\cdots\wedge y^{m}=\det(x^{i}\cdot y^{j})_{J},$$ то это
означает, что мы индуцируем из $\mathbb{E}$ в $\Lambda^{m}
(\mathbb{E})$ скалярное произведение поливекторов. Если взять
квадрат поливектора $x^{1}\wedge\cdots\wedge x^{m}$, т. е. его
скалярное произведение на себя, тогда мы получим определитель
Грама системы векторов $(x^{1},\ldots,x^{m})$, равенство нулю
которого эквивалентно линейной зависимости этой системы.
Линеаризация определителя Грама, т. е. корень квадратный из его
абсолютного значения, представляет собой полилинейную (но не
знакопеременную) функцию, которую обычно отождествляют с объемом
параллелепипеда, натянутого на систему $(x^{1},\ldots,x^{m})$.
Напомним, что скалярное произведение устанавливает изоморфизм
между $\mathbb{E}$ и дуальным к нему пространством
$\mathbb{E}^{*}$ с базисом $(\varepsilon^{i})_{I}$, двойственным к
ортонормированному базису $(e^{i})_{I}$, т. е.
$$\ast:\mathbb{E}\rightarrow\mathbb{E}^{*}:x\rightarrow
x^{*}:\langle x^{*},y\rangle=x\cdot y\quad(\forall\:x,y\in
\mathbb{E}):$$ $$\langle e_{i}^{*},e_{i}\rangle=e_{i}\cdot
e_{i}:e_{i}^{*}=e_{i}\cdot e_{i}\varepsilon_{i},$$ который
индуцирует изоморфизм
$\mathbb{T}^{m}(\mathbb{E})\simeq\mathbb{T}^{m}(\mathbb{E^{*}})$.
В частности, изоморфизм
$\Lambda^{m}(\mathbb{E})\simeq\Lambda^{m}(\mathbb{E^{*}})$
устанавливается соответствием между m--векторами и m--формами по
формуле $$\ast:\sum x_{J}e^{J}\rightarrow\sum
x_{J}^{*}\varepsilon^{J}:\sum_{\{(J)\}}x_{(J)}
e^{(i_{1}}\wedge\cdots\wedge
e^{i_{m})}\rightarrow\sum_{\{(J)\}}x_{(J)}
e^{(*i_{1}}\wedge\cdots\wedge e^{*i_{m})}.$$ Напомним также, что
евклидова структура линейного пространства задает оператор
дуализации Ходжа
$$\star:\Lambda^{m}(\mathbb{E})\rightarrow
\Lambda^{n-m}(\mathbb{E}):$$ $$\sum x_{J}e^{J}\rightarrow \langle
*\sum x_{J}e^{J},\epsilon\rangle: \sum x_{J}e^{J}\rightarrow \sgn
/J,J'/ \sum x_{J}^{*}e^{J'},$$ где дискриминант
$\epsilon:=e^{1}\wedge\cdots\wedge e^{n}$ при сворачивании с
m--формой $e^{J}$ приобретает знак подстановки $e^{J} \wedge
e^{J'}:=e^{(i_{1}}\wedge\cdots\wedge e^{i_{m})} \wedge
e^{(i_{m+1}}\wedge\cdots\wedge e^{i_{n})}$. Впрочем, в евклидовом
пространстве с положительной сигнатурой, где $x_{J}^{*}=x_{J}$, мы
получим $\star:\Lambda^{m}(\mathbb{E})\rightarrow
\Lambda^{n-m}(\mathbb{E}):\sum x_{J}e^{J}\rightarrow\sgn /J,J'/
\sum x_{J}e^{J'}$. Тем самым, наше определение векторного
произведения через четность частично упорядоченных подстановок
согласовано с определением оператора дуализации в евклидовом
пространстве с положительной сигнатурой через сворачивание с
дискриминантом.

Обратимся теперь к гладким дифференциальным формам. Пусть дано
поле гладких дифференциальных m--форм $a(x)=\sum a_{J}(x)dx^{J}$,
где $J:=(J),\quad dx^{J}:=dx^{(i_{1}}\wedge\cdots\wedge
dx^{i_{m})}$ и $a_{J}(x)$ --- произвольные гладкие отображения из
$\mathbb{L}$ в $\mathbb{R}$. Напомним, что внешним произведением
дифференциальной m--формы $a(x)=\sum a_{J}(x)dx^{J}$ и
дифференциальной 1--формы $b(x)=\sum_{n}b_{i}(x)dx^{i}$ называют
(m+1)--форму $a(x)\wedge b(x)=\sum b_{i}(x)a_{J}(x) \sgn \check{J}
dx^{\check{J}}$, где $\check{J}$ есть упорядоченная подстановка
присоединения $/(J)i/$, т.е. внешнее произведение получается в
резултате внешнего произведения этих форм и подстановок
присоединения, преобразующих базисные элементы $dx^{J}\wedge
dx^{i}$ к упорядоченному виду $dx^{\check{J}}$.

Вместе с тем, помимо внешнего произведения, существует внутреннее
произведение дифференциальной m--формы $a(x)=\sum a_{J}(x)dx^{J}$
на векторное поле $\bar{b}(x)=\sum_{n} b_{i}\frac{\partial
}{\partial x_{i}}$, под которым понимают свертку поля формы с
полем вектора. Однако, несмотря на простоту такого определения
внутреннего произведения, обратим ваше внимание на то, что для
сворачивания необходимо, выполнив подстановки изъятия, привести
форму $a(x)=\sum a_{J}(x)dx^{J}$ к виду $\sum a_{J}(x) \sgn
\hat{J} dx^{\hat{J}}\wedge dx^{i}$, а затем, заменив в этом
представлении внешнее произведение на $dx^{i}$ на произведение на
число $b_{i}(x)$, получить уже искомое внутреннее произведение
$\langle a(x),\bar{b}(x)\rangle= \sum b_{i}(x) a_{J}(x) \sgn
\hat{J} dx^{\hat{J}}$.

Внешняя производная гладкой дифференциальной формы (дифференциал)
задается формулой $da(x)= \sum d_{\check{J}}a(x)dx^{\check{J}}
=\sum da_{J}(x)\wedge dx^{J}$, где $da_{J}(x)=\sum_{n}
\frac{\partial a_{J}(x)}{\partial x_{i}}$. Тогда, в силу свойства
антикоммутативности подстановок упорядочивания с присоединением,
которые возникают при повторном дифференцировании, мы получим
тождество $dda(x)=0$, характеризующее внешнюю производную как
повышающий граничный гомоморфизм.

Установим теперь еще одно важное свойство внешеней производной.
Пусть дан такой произвольный предельно малый вектор $\Delta x =
\sum_{n}\Delta x^{i}=\sum_{n}\Delta x_{i}e^{i}$, что из него можно
сформировать ненулевой (m+1)--вектор $\Delta=\sum \Delta
x^{\check{J}}$, где $\Delta x^{\check{J}}: = \Delta
x^{(i_{1}}\wedge\cdots\wedge\Delta x^{i_{m+1})}$, и
соответствующий ему m--вектор $\delta\Delta =\sum\Delta x^{J}$,
где $\Delta x^{J}=\Delta x^{(i_{1}}\wedge\cdots\wedge\Delta
x^{i_{m})}$. Тогда $\langle da(x),\Delta \rangle= \sum\langle
d_{\check{J}}a(x)dx^{\check{J}},\Delta x^{\check{J}} \rangle$, где
$\langle d_{\check{J}}a(x)dx^{\check{J}},\Delta x^{\check{J}}
\rangle =d_{\check{J}}a(x)\cdot\Delta x_{i_{1}}\cdots\Delta
x_{i_{m+1}}$. С другой стороны, разложив в этом выражении
компоненту дифференциала и компоненту поливектора, мы будем иметь
выражение $\langle d_{\check{J}}a(x)dx^{\check{J}},\Delta
x^{\check{J}} \rangle =\langle \sum_{n}\frac{\partial
a_{J}(x)}{\partial x_{i}}\wedge dx^{J},\sum_{n}\Delta x^{i} \wedge
\Delta x^{J} \rangle$, и если затем свернем 1--форму и вектор, то
получим равенство $\langle d_{\check{J}}a(x)dx^{\check{J}}, \Delta
x^{\check{J}} \rangle = \langle \sum_{m+1}[a_{J}(x+\Delta x^{i}) -
a_{J}(x) )]dx^{J}, \sum_{m+1}\Delta x^{J} \rangle$. Таким образом,
доказано утверждение
\begin{Proposition}
$\langle da(x),\Delta \rangle=\langle [a(x+\Delta x)-a(x)],
\delta\Delta \rangle$
\end{Proposition}

Данное предложение мы применим в следующем разделе при выводе
формулы Стокса, где предельно малые вектора служат остовом для
предельно малого параллелепипеда а интегрирование по поверхности
отождествляется с суммированием сверток диференциальной формы и
поливекторов, соответствующих этим параллелепипедам. Сейчас же,
заметим только, что величина $\langle da(x),\Delta \rangle$
характеризует линейную часть приращения свертки дифференциальной
формы в точке $x$ и поливектора, соответствующего граням
параллелепипедов, построенных на предельно малых векторах, при
сдвиге каждой грани на противоположную и соответствующем изменении
значения дифференциальной формы.

\section{Параллелепипеды}

Начнем мы с классического определения параллелепипедов, которое
здесь повторим, сославшись на университетский учебник \cite{вин}.
\begin{defn}
Параллелепипедом $\pi^{m}$, натянутым на систему линейно
независимых векторов $x^{1},\ldots ,x^{m}$, выходящих из точки
$x^{0}$ n--мерного аффинного пространства $\mathbb{R}^{n}$,
называется множество точек $x^{0}+\sum_{m}\alpha_{i}x^{i}$,
задаваемое условием $0\leq\alpha_{i}\leq 1\quad\forall\:i\in J$.
\end{defn}
Затем, в дополнение к этому геометрическому понятию, мы дадим
алгебраическое понятие ориентированного параллелепипеда.
\begin{defn}
Ориентированным параллелепипедом называется параллелепипед, в
котором определенным образом задано направление всех цепочек ребер
$x^{/i_{1}},\ldots,x^{i_{m}/}$. Направление цепочки от $x^{0}$ к
$x^{0}+\sum_{m}x^{i}$ считается положительным, а обратное ---
отрицательным. Если в параллелепипеде всем четным цепочкам
(цепочкам с четными последовательностями ребер) задать
положительное направление, а нечетным
--- отрицательное, то такой параллелепипед называется положительно
ориентированным и обозначается $+\pi^{m}$, т. е. $\sgn\pi^{m}=+1$.
В обратном случае параллелепипед называется отрицательно
ориентированным и обозначается $-\pi^{m}$, т. е. $\sgn\pi^{m}=-1$.
\end{defn}
На основании геометрического определения мы полагаем, что гранями
параллелепипеда $\pi^{m}$ служат $2m$ параллелепипедов
$\pi^{m-1}$, а именно, $m$ параллелепипедов $\pi^{m-1}_{*j}$,
натянутых на подсистемы $m-1$ векторов, полученных исключением
вектора $x^{j}$, где $j\in J$, из системы $x^{1},\ldots ,x^{m}$, и
выходящих из начальной точки $x^{0}$, а также $m$ параллелепипедов
$\pi^{m-1}_{j*}$, выходящих из $m$ точек $x^{0}+x^{i}$ и натянутых
соответственно на подсистемы без $x^{j}$ каждый. Таким образом, у
параллелепипеда $2^{m}$ вершин, $2m$ граней, и $m2^{m-1}$ ребер.
Заметим также, что все грани граней параллелепипеда $\pi^{m}$
соединяют соответствующие пары соседних граней, т. е.
$\pi^{m-1}_{j*}\cap\pi^{m-1}_{k*}=\pi^{m-2}_{jk*},\quad
\pi^{m-1}_{*j}\cap \pi^{m-1}_{*k}=\pi^{m-2}_{*jk},\quad
\pi^{m-1}_{j*}\cap\pi^{m-1}_{*k}=\pi^{m-2}_{j*k}$, где $j\neq k$.

\begin{defn}
Ориентированной поверхностью параллелепипеда называется
совокупность определенным образом ориентированных его граней.
Ориентация граней $\pi^{m-1}_{*j}$ параллелепипеда соответствует
четности подстановок $/(1,\ldots,m)\hat{j}/$ а грани
$\pi^{m-1}_{j*}$, противоположные граням $\pi^{m-1}_{*j}$, имеют
также и противоположную ориентацию. Иначе говоря,
$\sgn\pi^{m-1}_{*j}
=-\sgn\pi^{m-1}_{j*}=\sgn/(1,\ldots,m)\hat{j}/$.
\end{defn}
Тогда, в силу антикоммутативности подстановок с изъятием, мы
получим одно замечательное свойство

\begin{Proposition}\label{grani}

Общие грани соседних граней ориентированной поверхности
параллелепипеда имеют в них противоположную ориентацию.

\end{Proposition}

Действительно, в параллелепипеде $\pi^{m-1}_{*j}$ ориентация грани
$\pi^{m-2}_{*jk}$ равна $\sgn/(1,\ldots,m)\hat{j}\hat{k}/$ а в
параллелепипеде $\pi^{m-1}_{*k}$ ориентация этой же грани равна
$\sgn/(1,\ldots,m)\hat{k}\hat{j}/$. Аналогично, в параллелепипеде
$\pi^{m-1}_{j*}$ ориентация грани $\pi^{m-2}_{jk*}$ равна
$\sgn/(1,\ldots,m)\hat{j}\hat{k}/$ а в параллелепипеде
$\pi^{m-1}_{k*}$ ее ориентация равна $\sgn/(1,\ldots,m)
\hat{k}\hat{j}/$. Вместе с тем, в параллелепипеде $\pi^{m-1}_{j*}$
ориентация грани $\pi^{m-2}_{j*k}$ равна $-\sgn/(1,\ldots,m)
\hat{j}\hat{k}/$ а в параллелепипеде $\pi^{m-1}_{*k}$ ее
ориентация равна $-\sgn/(1,\ldots,m)\hat{k}\hat{j}/$.

Алгебраическая конструкция ориентировнного параллелепипеда
допускает формализацию. Действительно, если допустить кратность
ориентированного параллелепипеда, реализуемую посредством
расслоения цепочек его ребер, то алгебраическую сумму чисел
положительно и отрицательно ориентированных слоев можно задавать
некоторым произвольным целым числом. Тем самым, мы можем говорить
о z-ориентированном параллелепипеде $z\pi^{m}$, который является
элементом абелевой группы $\Pi_{m}$, изоморфной $\mathbb{Z}$.

В свою очередь, геометрическая конструкция параллелепипеда через
возможность формирования гомологической сети параллелепипедов
допускает расширение до понятия клеточной поверхности.
Действительно, параллелепипед $\pi^{m}$ в пространстве
$\mathbb{R}^{n}$ может быть продолжен геометрически, через
собственные $(m-1)$-мерные грани так, что новый $m$-мерный
параллелепипед натягивается на произвольную, но уже существующую
грань, и некоторый новый вектор, не лежащий в подпространстве этой
грани. Если допустить, что параллелепипеды гомеоморфно
отображаются в клетки, которые сами себя не пересекают, тогда
конечное или счетное объединение клеток-параллелепипедов
$\pi^{m}$, произвольные пары которых либо вовсе не пересекаются,
либо пересекаются своими общими гранями произвольной размерности,
можно отождествить с $m$-мерным клеточным комплексом. Односвязной
$m$-мерной клеточной поверхностью $\Omega^{m}$ мы назовем такой
$m$-мерный клеточный комплекс, в котором каждая пара
клеток-параллелепипедов может быть соединена в цепь, связанную
общими $(m-1)$-мерными гранями ее звеньев -- соседних
клеток-параллелепипедов. Если все клетки-параллелепипеды
односвязной клеточной поверхности пересекаются только своими
общими гранями произвольной размерности и каждая общая грань
размерности $(m-1)$ является пересечением не более 2
клеток-параллелепипедов, то мы говорим об односвязной клеточной
поверхности без самопересечений или о клеточном многообразии.
Границей $d\Omega^{m}$ клеточногй поверхности $\Omega^{m}$ мы
назовем объединение всех $(m-1)$-мерных граней
клеток-параллелепипедов поверхности $\Omega^{m}$ за вычетом ее
общих граней. Поскольку все грани граней параллелепипеда являются
общими гранями, то и все грани (m-1)-мерных параллелепипедов
границы являются общими гранями, а следовательно граница является
объединением односвязных (m-1)-поверхностей без границ.

Топологически инвариантным (гомеоморфным) преобразованием
клеточной поверхности с границей мы назовем такое присоеденение к
ее границе новой клетки-параллелепипеда или изъятие уже
существующей клетки-параллелепипеда, которое не нарушает связности
примыкающей к этой клетке границы, т.е. не изменяет числа
односвязных кусков ее границы.

Множество всех m--мерных клеток-параллелепипедов, составляющих
поверхность $\Omega^{m}$, порождает свободную абелеву группу
$\Pi_{m}(\Omega^{m})$, образованную формальными суммами
$\sum_{K}z_{i}\pi^{m}_{i}$, где $z_{i}\in \mathbb{Z},\quad
K\subseteq\mathbb{N}$. Сумму, составленную из $\pm
1$--ориентированных клеток-параллелепипедов поверхности
$\Omega^{m}$, грани пересечения соседних клеток-параллелепипедов
которой имеют в них противоположную ориентацию, мы назовем
ориентированной клеточной поверхностью. Поверхность $\Omega^{m}$,
которая допускает ориентированную поверхность, называется
ориентируемой поверхностью. Заметим при этом, что при любом
гомеоморфном присоединении или изъятии клетки-параллелепипеда
грань соединения примыкающих клеток всегда имеет в них
противоположную ориентацию. Действительно, если совместить
примыкающие клетки сохраняющим ориентацию сдвигом, то грань
примыкания отобразится в две противоположные грани, имеющие
противоположную ориентацию. Следовательно ориентация гомеоморфно
присоединяемой (изымаемой) клетки согласована с ориентацией
соседних клеток, и поэтому мы получим

\begin{Proposition}
Ориентируемость является топологическим инвариантом поверхности
параллелепипедов $\Omega^{m}$
\end{Proposition}

Кроме того, всякая клеточная поверхность $\Omega^{m}$ порождает
свободную абелеву группу $\Pi_{m-1}(\Omega^{m})$, образованную
формальными суммами граней ее клеток-параллелепипедов. Группу
$\Pi_{m-1}(\Omega^{m})$ мы называем группой граней поверхности.
Группа граней одной клетки-параллелепипеда $\Pi_{m-1}(\pi^{m})$
состоит из формальных сумм $\sum_{2m}z_{i}\pi^{m-1}_{i}$ а
отображение $$\gamma:\pi^{m}\rightarrow \Pi_{m-1}(\pi^{m}):$$
$$\pi^{m}\rightarrow \sum_{m}\sgn/(1,\ldots,m)\hat{j}/\cdot
\pi^{m-1}_{*j}-\sum_{m}\sgn/(1,\ldots,m)\hat{j}/\cdot\pi^{m-1}_{j*}$$
индуцирует граничный гомоморфизм соответствующих абелевых групп
$$\delta:\Pi_{m}(\Omega^{m}) \rightarrow \Pi_{m-1}(\Omega^{m}):
\sum_{K} z_{i}\pi^{m}_{i} \rightarrow \sum_{K}z_{i}
\gamma(\pi^{m}_{i}).$$ Действительно, в силу предложения
\ref{grani} имеет место тождество $\gamma\gamma(\pi^{m})=0$,
применяя которое покомпонентно, получим другое тождество
$\delta\delta(\Omega^{m})=0$, определяющее граничное отображение и
цепной комплекс $\{\Pi_{m-1}(\Omega^{m}),\delta\}$.

Положим $$Z_{m-1}(\Omega^{m}):=\Ker\{\delta:\Pi_{m-1}(\Omega^{m})
\rightarrow\Pi_{m-2}(\Omega^{m})\}\subset\Pi_{m-1}(\Omega^{m})$$ и
назовем $Z_{m-1}(\Omega^{m})$ группой циклов. Положим также
$$B_{m-1}(\Omega^{m}):=\im\{\delta:\Pi_{m}(\Omega^{m})
\rightarrow\Pi_{m-1}(\Omega^{m})\}\subset\Pi_{m-1}(\Omega^{m})$$ и
назовем $B_{m-1}(\Omega^{m})$ группой границ. Так как
$\delta\delta=0$, то $B_{m-1}(\Omega^{m})\subset
Z_{m-1}(\Omega^{m})$. Факторгруппу $$H_{m-1}(\Omega^{m}):=
Z_{m-1}(\Omega^{m})/B_{m-1}(\Omega^{m})$$ назовем $(m-1)$-мерной
группой гомологий. Тогда, поскольку гомеоморфное присоединение или
изъятие клетки-параллелепипеда добавляет в группу граней или
удаляет из нее элементарный цикл, который является границей, то мы
имеем

\begin{Proposition}
Группа $(m-1)$-мерных гомологий является топологическим
инвариантом клеточной поверхности $\Omega^{m}$
\end{Proposition}

Таким образом, применение понятия ориентированного параллелепипеда
позволяет сделать естественный переход к построению стандартной
теории гомологий. Однако уже в следующем разделе мы покажем, что
концепция применения параллелепипедов с успехом может работать и в
нестандартной топологической теории.

Обратимся теперь к вопросам измерения, связанным с поверхностями.
Поскольку внешнее произведение векторов, на которые натягивается
параллелепипед, задает отображение
$\bar{\pi}^{m}:\pi^{m}\rightarrow\wedge^{m}(\mathbb{L})$,
индуцирующее эквивалентность элементов пространства $\{\pi^{m}\}$,
то прострвнство $\wedge^{m}(\mathbb{L})$ и
$\wedge^{m}(\mathbb{L^{*}})$ можно интерпретировать как
пространство эквивалентных приборов и воздействий на них
соответственно. Иначе говоря, мы полагаем, что свертка
$\langle\bar{x},\tilde{y}\rangle$, где $\bar{x} \in
\wedge^{m}(\mathbb{L}),\quad\tilde{y}\in
\wedge^{m}(\mathbb{L^{*}})$, есть результат измерения воздействия
$\tilde{y}$ на прибор $\bar{x}$. В соответствии с такой
концепцией, всякой поверхности $\Omega^{m}$ и полю
дифференциальных m--форм $a(x)$ можно сопоставить сумму сверток
$\sum_{K}\langle a(x^{0}_{i}),\bar{\pi}^{m}_{i}\rangle$, которую
следует понимать как результат воздействия на поверхность,
составленную из $K$ приборов--параллелепипедов. В том случае,
когда поверхность $\Omega^{m}$ составлена из предельно малых
параллелепипедов $\Delta\pi^{m}$, мы назовем ее интегральной
поверхностью а сумму сверток обозначим символом интеграла
$\int_{\Omega^{m}}\langle a(x),\Delta\bar{\pi}^{m}\rangle$.
Заметим при этом, что для канонической поверхности, все
параллелепипеды которой построены на векторах, колинеарных
базисным, принято обозначение $\int_{\Omega^{m}}a(x)dx^{J}$.
Посколку общие грани являются одновременно начальными и конечными
порождающими гранями в соседних параллелепипедах, то применив
разложение
$$\langle da(x),\Delta\bar{\pi}^{m}\rangle=\sum_{m}\langle
a(x+\Delta x^{j}),
\Delta\bar{\pi}^{m-1}_{j*}\rangle-\sum_{m}\langle a(x),
\Delta\bar{\pi}^{m-1}_{*j}\rangle$$ ко всем предельно малым
параллелепипедам интегральной поверхности, мы получим формулу
Стокса $$\int_{\Omega^{m}}\langle da(x),\Delta\bar{\pi}^{m}
\rangle= \int_{\delta\Omega^{m}}\langle
a(x),\Delta\bar{\pi}^{m-1}\rangle.$$

Если отображение $\bar{\pi}^{m}$ индуцировано векторным
произведением, то с ориентированной поверхностью $\Omega^{m}$
можно связать векторный объем $\vol(\Omega^{m}):=
\sum_{K}\bar{\pi}^{m}_{i}$. Наконец, со всякой интегральной
поверхностью, заданной в евклидовом пространстве положительной
сигнатуры, можно связать линейный объем
$\int_{\Omega^{m}}\sqrt{\abs{\Delta\bar{\pi}^{m}\cdot
\Delta\bar{\pi}^{m}}}$, который в наблюдаемом пространстве
ассоцируется с длиной, плошадью и собственно объемом.

\section{Ориентируемые многообразия}

Прежде всего дадим ряд определений. Пусть неориентированный граф
без петель так вложен в пространство $\mathbb{R}^{n+1}$, что в
каждой из вершин собран пучек из $n$ криволинейных ребер, при этом
все ребра разбиты на $n$ классов а в каждой вершине встречается по
одному представителю каждого класса, причем все касательные
вектора к ребрам в каждой вершине линейно независимы. Назовем
такой граф (вместе с его укладкой в $\mathbb{R}^{n+1}$) каркасом и
натянем на него поверхность с границей, т.е. образуем клеточную
$n$-поверхность с каркасом в качестве одномерного остова
$n$-мерного клеточного комплекса. Определим же мы наш $n$-мерный
клеточный комплекс по индукции: 1-мерный клеточный комплекс задан
своим каркасом, 2-мерный клеточный комплекс задается гомеоморфными
параллелограму 2-мерными клетками с границами, состоящими из
всевозможных ребер двух разных классов, $n$-мерный комплекс
задается гомеоморфной $n$-мерному параллелепипеду клеткой с
границей, состоящей из $(n-1)$-мерного комплекса. Натянутую на
каркас поверхность будем также называть поверхностью каркаса.
Тогда имеет место
\begin{Lemma}
Поверхность каркаса гомеоморфна $n$--мерному симплексу.
\end{Lemma}
\begin{proof}
С одной стороны, поскольку граф у нас связный, то соответствующий
ему каркас можно гомеоморфно стянуть в точку. С другой стороны,
любая вершина поверхности каркаса одновременно является вершиной
предельно малого элемента этой поверхности; т.е., гомеоморфного
точке $n$--мерного симплекса. Следовательно, стягивая каркас к
некоторой вершине, можно совместить его с этим предельно малым
симплексом.
\end{proof}

Вложим теперь в пространство $\mathbb{R}^{n+1}$ связку
окружностей, целиком лежащую в $n$ не параллельных друг другу
плоскостях, причем все узлы связки лежат на одной прямой а каждая
окружность связки пересекается в двух своих противоположных
точках, которые являются узлами пучков из $n$ окружностей, лежащих
в $n$ плоскостях. Назовем такую связку окружностей $c$--каркасом и
выделим в нем определенное семейство обычных каркасов с вершинами
в узлах $c$--каркаса и с ребрами, состоящими из полуокружностей
$c$--каркаса. Для этого мы каждый класс окружностей $c$--каркаса,
принадлежащих одной плоскости, разделим посредством прямой,
проходящей через его узлы, на два класса полуокружностей. Тогда,
выбирая из полуокружностей $c$--каркаса ребра в соответствии с
определенным выбором классов полуокружностей, мы получим семейство
из $2^{n}$ различных каркасов, которые будем называть
$s$--каркасами. В свою очередь, суммарную поверхность всех
$s$--каркасов мы назовем поверхностью $c$--каркаса. Тогда имеет
место
\begin{Lemma}\label{ostov}
Поверхность $c$--каркаса является $n$--мерным клеточным замкнутым
ориентируемым многообразием.
\end{Lemma}
\begin{proof}
Утверждение о том, что поверхность $c$--каркаса есть замкнутое
многообразие следует из того, что всякая часть границы поверхности
произвольного $s$--каркаса является одновременно частью границы
поверхности другого $s$--каркаса. Действительно, для поверхностей
двух $s$--каркасов, отличающихся только одним классом
полуокружностей, часть границы, полученная исключением этого
класса полуокружностей, является их общей границей. Тем самым,
поверхность $c$--каркаса склеена из поверхностей $s$--каркасов
(гомеоморфных шару) по их границам, но сама не имеет границы.
Ориентируемость нашего многообразия следует из того, что общая
часть границы поверхностей двух $s$--каркасов имеет в них
противоположную ориентацию, которая задается противоположным
классом полуокружностей, отличающим поверхности двух граничащих
друг с другом $s$--каркасов.
\end{proof}
Рассмотрим теперь вопрос о классификации $c$--каркасов. Если мы
обратимся к нашему определению $c$--каркасов, то увидим, что
классификация $c$--каркасов эквивалентна классификации
неориентированных графов без петель, у которых все ребра разбтиы
на $n$ классов а каждой вершине инцендентны $n$ смежных ребер,
принадлежащих различным классам. Действительно, достаточно
установить соответствие между ребрами одного класса и
окружностями, лежащими в одной плоскости, а также между вершинами
графа и узлами $c$--каркаса. Далее мы сформируем элементарный
$c$--каркас, соответствующий простейшему графу. Для этого,
произвольным образом разложим $n$ в сумму $n_{1}+\cdots +n_{m}=n$,
где $1\leq m \leq n$, и построим $2^{m}$--вершинный граф,
изоморфный $m$--параллелепипеду с ребрами кратности $n_{i}$, т. е.
такой граф, вершины которого совпадают с вершинами
$m$--параллелепипеда а всякие $n_{i}$ его кратных ребер
соответствуют каждому ребру параллелепипеда из соответствующего
класса его параллельных ребер. Тогда, поскольку у
$m$--параллелепипеда $m2^{m-1}$ ребер, то с учетом их кратности мы
получим $n2^{m-1}$ ребер графа или $n2^{m-1}$ окружностей
соответствующего $c$--каркаса. Элементарный $c$--каркас,
соответствующий разложению $n_{1},\ldots, n_{m}$ мы обозначим
символом $S^{n_{1}}\times\cdots \times S^{n_{m}}$ и заметим, что
из двух произвольных элементарных $c$--каркасов можно составить
новый $c$--каркас $S^{n_{1}}\times\cdots \times S^{n_{m}}+
S^{n_{1}}\times\cdots \times S^{n_{l}}$, которому соответствует
граф, полученный из двух простейших графов путем разрыва ребер в
одной произвольной вершине каждого из них и последующего
совмещения соответствующих ребер двух этих графов. Однако, если
$l=1$, то мы получим исключение, так как в этом случае
$S^{n_{1}}\times\cdots \times S^{n_{m}}+S^{n}=S^{n_{1}}
\times\cdots \times S^{n_{m}}$. Таким образом, мы получили полную
классификацию $c$--каркасов, которая полностью определяется
свободной абелевой группой с элементарными $c$--каркасами (кроме
$S^{n}$) в качестве образующих и каркасом $S^{n}$ в качестве
нулевого элемента.

Сейчас мы приступим к доказательству обратного \ref{ostov}
утверждения.
\begin{Theorem}
Всякое $n$-мерное замкнутое ориентируемое многообразие, вложенное
в $\mathbb{R}^{n+1}$, гомеоморфно поверхности некоторого
$c$--каркаса
\end{Theorem}
\begin{proof}
В двумерном случае это действительно так, поскольку классификация
$c$--каркасов совпадает с классификацией двумерных ориентируемых
замкнутых многообразий. В общем случае, ввиду отсутствия полной
классификации многообразий, нам следует установить, что любое
замкнутое ориентируемое многообразие, вложенное в
$\mathbb{R}^{n+1}$, обладает каркасом, гомеоморфным некоторому
$c$--каркасу. Пусть дано многообразие $M^{n}$ и его клеточное
покрытие, образующее клеточную поверхность. Возьмем произвольную
точку $x$ многообразия $M^{n}$, которая совпадает с вершиной
произвольной клетки покрытия, и заметим, что она является общей
вершиной для $2^{n}$ клеток-параллелепипедов, лежащих в ее
окрестности, а любое ребро, инцендентное этой вершине, является
общим ребром для $2^{n-1}$ клеток-параллелепипедов. Затем,
используя $2^{n}$ цветов краски, раскрасим покрытие многообразия
$M^{n}$ так, чтобы краски не перемешивались а первые мазки краски
покрыли бы каждую клетку окрестности точки $x$ своим цветом. Иначе
говоря, требуется, используя $2^{n}$ клеток в качестве начальных
для $2^{n}$ клеточных поверхностей, гомеоморфным присоединением
новых клеток разбить все покрытие многообразия $M^{n}$ на $2^{n}$
клеточных поверхностей. Раскрасив таким образом покрытие
многообразия $M^{n}$, мы получим некий каркас. Действительно,
точки многообразия, в которых пересекаются поверхности всех
цветов, будут служить вершинами, а линии многообразия, в которых
пересекаются $2^{n-1}$ поверхностей различного цвета, будут
служить ребрами нашего каркаса. Но поскольку исходное многообразие
было ориентируемым, то полученный каркас допускает согласованную
ориентацию ребер в каждой вершине. Следовательно всякую пару
противоположных ребер этого каркаса, выходящие из одной вершины,
можно сопоставить окружности некоторого $c$--каркаса. Тем самым,
проблема полной классификации замкнутых ориентируемых
$n$-многообразий, вложенных в $\mathbb{R}^{n+1}$, находит решение,
соответствующее классификации $c$--каркасов.
\end{proof}
В заключение заметим, что многообразие, имеющее $c$--каркас типа
$S^{n}$, гомеоморфно сфере $S^{n}$, а многообразие, имеющее
$c$--каркас, представленный суммой двух элементарных
$c$--каркасов, образуется из элементарных многообразий с помощью
вырезания в каждом из них по клетке и склеивания их по границе
клетки. С учетом этих замечаний справедливо утверждение.
\begin{Proposition}
Фундаментальная группа замкнутого ориентируемого многообразия
тривиальна, если во всех его образующих элементарных каркасах типа
$S^{n_{1}}\times\cdots \times S^{n_{m}}$ не встречается ни одной
компоненты $S^{1}$.
\end{Proposition}
\begin{proof}
Действительно, при данных условиях исходное многообразие можно
склеить из элементарных многообразий, каждое из которых можно
разложить в произведение подмногообразий, гомеоморфных некоторым
сферам $S^{k}$, где каждое $k>1$. Но поскольку фундаменнтальная
группа произведения таких сфер тривиальна а каждая граница склейки
гомеоморфна сфере $S^{n-1}$, которая в случае $n>2$ имеет
тривиальную фундаментальную группу, то утверждение доказано.
\end{proof}
Из последнего предложения и из того, что классификация замкнутых
ориентируемых многообразий в размерности три сводится к
классификации соответствующих каркасов, которые могут быть
составлены из каркасов типа $S^{1}\times S^{1}\times S^{1}$ и
$S^{1}\times S^{2}$, где обязательно встречается компонента
$S^{1}$, либо это есть каркас типа $S^{3}$, следует, что замкнутые
ориентируемые 3-многообразия, имеющие тривиальную фундаментальную
группу, гомеоморфны 3-мерной сфере. Таким образом, классическая
гипотеза Пуанкаре, сформулированная в классе вложенных
многообразий, является простым следствием классификации замкнутых
ориентируемых многообразий посредством классификации их каркасов.

\appendix
\section{Минимальные потоки и поверхности}

Векторные поля и минимальные поверхности, заданные в евклидовом
пространстве, основательно описаны в учебной литературе,
например,\cite{AMI} \cite{гор}, \cite{бла} и \cite{тор}, сейчас же
мы добавим к известному лишь то новое, что вытекает из нашей
алгебры параллелепипедов, либо посмотрим на известные факты
несколько с другой стороны.

Пусть в $\mathbb{R}^{n}$ задано произвольное гладкое ковекторное
поле, т. е. поле дифференциальной 1--формы
$a(x)=\sum_{n}a_{j}(x)dx^{j}$. Тогда мы имеем также дуальную к
$a(x)$ дифференциальную (n-1)--форму $\star a(x)=\sgn J
\sum_{n}a_{j} (x)dx^{J}$, где $J:=/I,\hat{j}/$. Дифференциальные
формы $a(x)$ и $\star a(x)$ можно представить стационарным потоком
некой идеальной жидкости, измеряемым 1--мерными и (n-1)--мерными
приборами (поливекторами) в результате их сворачивания с $a(x)$ и
$\star a(x)$ соответственно. Если в качестве приборов для
измерения потока используются 1--мерные или (n-1)--мерные
интегральные поверхности, то в результате измерения мы получим
интегральные суммы $\int_{\Omega^{1}}\langle a(x),dx\rangle$ и
$\int_{\Omega^{n-1}}\langle \star a(x),dx^{n-1}\rangle$
соответственно, где $dx:=\Delta\bar{\pi}^{1}$ и
$dx^{n-1}:=\Delta\bar{\pi}^{n-1}$. Применение предельно малых
поверхностей без границ дает нам возможность установить локальные
характеристики стационарного потока. Действительно, отсутствие
локального вращения потока $\oint\langle a(x),dx^{1}\rangle=0$
обеспечено дифференциальным условием $da(x)=0$, а отсутствие
стоков и истоков $\oint\langle\star a(x),dx^{n-1}\rangle=0$
обеспечено условием $d\star a(x)=0$.

Критерии замкнутости потока мы дополним критериями его
устойчивости. Так, ламинарное течение стационарного потока
удовлетворяет дифференциальному условию голономности $a(x)\wedge
da(x)=0$. Заметим при этом, что всякая точная 1--форма
$a(x)=d\varphi(x)$ голономна а ее семейство интегральных
поверхностей задано поверхностями уровня функции $\varphi(x)$. В
свою очередь, всякое голономное ковекторное поле $a(x)$ колинеарно
некоторому потенциальному ковекторному полю $d\varphi(x)$, т.е.
$a(x)=k(x)d\varphi(x)$, где $k(x)$ --- произвольная гладкая
функция. Действительно, $da(x)= d[k(x)d\varphi(x)]=dk(x)\wedge
d\varphi(x)$, и поэтому $a(x)\wedge da(x)=k(x)d\varphi(x) \wedge
dk(x)\wedge d\varphi(x)=0$, обратно, если уравнение $a(x)=0$
вполне интегрируемо, то его интегральным поверхностям всегда можно
поставить в соответствие поверхности уровня некоторой функции.

Равновесное течение ламинарного стационарного потока удовлетворяет
соответствующему дифференциальному условию. Так, для
потенциального поля $a(x)=d\varphi(x)$ это будет уравнение
$\Delta\varphi(x)=0$. Действительно, величина $\Delta\varphi(x)=
d\star d\varphi(x)$ характеризует разность давлений потока на
предельно малые элементы поверхностей уровня, полученные
градиентным переносом в предельно малом элементе объема, а именно,
$$\int_{\Delta V}\langle d\star d\varphi(x), dx^{n} \rangle =
\int_{\Delta \Phi_{2}}\langle \star d\varphi(x), dx^{n-1} \rangle
- \int_{\Delta \Phi_{1}}\langle \star d\varphi(x), dx^{n-1}
\rangle,$$ где учтено отсутствие давления на боковую поверхность
предельно малого элемента объема $\Delta V$, образованного
заметанием предельно малой площадки $\Delta \Phi_{1}$ поверхности
уровня $\varphi(x)=c$ на площадку $\Delta \Phi_{2}$ поверхности
уровня $\varphi(x)=c+h$ при ее переносе градиентным векторным
полем. Следовательно гармоничность функции уровня свдетельствует
об отсутствии разности давлений, т. е. о равновесии (устойчивости)
соответствующего стационарного потока. В свою очередь, если
$a(x)=k(x)d\varphi(x)$, то в качестве критерия устойчивости потока
следует рассматривать уравнение $d\star k(x)d\varphi(x)=0$.
Заметим при этом, что поскольку $d\star da(x)= d\star dk(x)\wedge
\star d\varphi(x)\pm \star dk(x)\wedge d\star d\varphi(x)$, то
равенство нулю этого выражения равнозначно гармоничности функций
$\varphi(x)$ и $k(x)$.

Покажем теперь, что гармонические функции, порождающие
потенциальные ковекторные поля $a(x)=d\varphi(x)$, обеспечивают
локальный экстремум криволинейного интеграла $\int\langle
a(x),dx\rangle$, определенного на произвольной кривой с
фиксированными границами, и в этом смысле определяют минимальные
стационарные потоки. Действительно,
$$\Lvar\int\langle a(x),dx\rangle= \delta \int_{x}^{x+\Delta
x}\langle a(x),dx\rangle= \delta \langle a(x+\Delta x)-a(x),\Delta
x \rangle=0,$$ и поэтому интегральная вариационная задача сводится
к вариационной задаче на минимум приращения свертки. Вместе с тем,
поскольку $a(x+\Delta x)=a(x)+\sum_{n}\langle da_{j}(x),\Delta
x\rangle dx^{j}$, то приращение свертки равно квадратичной форме,
т. е. $\langle a(x+\Delta x)-a(x),\Delta x \rangle=(H\Delta
x,\Delta x)$, где оператор $H$ в матричном представлении задается
матрицей Гессе $\left(\frac{\partial^{2}\varphi(x)}{\partial
x_{i}\partial x_{j}}\right)$. Таким образом, вариационное
уравнение $\delta (H\Delta x,\Delta x)=0$ сводит нашу задачу к
поиску минимального самосопряженного оператора. В этой связи
напомним \cite{гел}, что всякий самосопряженный оператор $C$ с
точностью до отражений может быть разложен на оператор $A$ с
нулевым следом и положительно определенный оператор $B$, т. е.
$C=\pm(A+B)$. Действительно, достаточно привести матрицу оператора
$C$ к диагональному виду и разложить ее диагональные элементы. В
то же время, для квадратичных форм выполняется неравенство
$(Ax,x)<((A+B)x,x)$, следовательно самосопряженный оператор с
нулевым следом будет минимальным оператором и решение $\tr
H=\Delta\varphi(x)\equiv 0$, определяет минимальное ковекторное
поле как градиент гармонической функции. Заметим также, что такое
же решение имеет вариационная задача на экстремум интеграла
$\int\langle \star a(x),dx^{n-1}\rangle$ по произвольной
(n-1)--поверхности с фиксированными границами при локальной
вариации точной 1--формы. Действительно, вариационное уравнение
$$\Lvar\int\langle \star a(x),dx^{n-1}\rangle= \delta \langle\star
a(x+\Delta x)-\star a(x),\star\Delta x \rangle=0$$ не меняет
постановки вариационной задачи, так как свертка вектора с
ковектором равна свертке их дуализаций.

Пусть теперь в евклидовом пространстве положительной сигнатуры,
где мы не делаем отличия векторов от ковекторов, дано поле точной
1--формы $d\varphi(x)$. Рассмотрим производное от него единичное
голономное поле $n(x)= k(x)d\varphi(x)$, где
$k(x)=1/|d\varphi(x)|$, и вычислим площадь $S$ поверхности
$\Phi^{n-1}$, которая является ограниченной частью поверхности
уровня функции $\varphi(x)$. Покажем, что
$$S(\Phi^{n-1})=\int_{\Phi^{n-1}}(\star n(x),dx^{n-1}).$$
Действительно, единичный вектор $n(x)$ ортогонален к плоскости,
касательной к поверхности уровня функции $\varphi(x)$ в точке $x$
а единичный поливектор $\star n(x)$ коллинеарен предельно малому
поливектору $\Delta\bar{\pi}^{n-1}(x)$, и поэтому площадь
последнего равна их скалярному произведению $(\star n(x),\Delta
\bar{\pi}^{n-1}(x))$, откуда следует, что
$$S(\Phi^{n-1})=\sum_{\Phi^{n-1}}S(\Delta
\pi^{n-1})=\int_{\Phi^{n-1}}(\star n(x),dx^{n-1}).$$ Заметим
также, что
$$\int_{\Delta V}\langle d\star n(x), dx^{n} \rangle =
S(\Delta\Phi_{2})-S(\Delta\Phi_{1}),$$ где площадка
$\Delta\Phi_{2}$ образована ортогональным переносом площадки
$\Delta\Phi_{1}$ на предельно малое расстояние с поверхности
уровня $\Phi_{1}$ на поверхность уровня $\Phi_{2}$, а
следовательно
$$\lim_{\Delta V\rightarrow 0}\frac{S(\Delta\Phi_{2})-
S(\Delta\Phi_{1})}{\Delta V}=d\star n(x).$$

Вместе с тем, легко показать, что локальное решение вариационного
уравнения
$$\Lvar\int_{\Phi^{n-1}}(\star n(x),dx^{n-1})=0$$
приводит к дифференциальному уравнению минимальности $d\star
n(x)=0$. Тем самым, локально минимальные параметризованные
поверхности евклидова пространства с нулевой средней кривизной и
поверхности уровня функций, удовлетворяющих уравнению
минимальности, эквивалентны. Впрочем, эта эквивалентность следует
также из того хорошо известного факта, что средняя кривизна
поверхности уровня, заданной в n--мерном евклидовом пространстве,
вычисляется по формуле $H(x)=
-(\frac{1}{n})d\star(\frac{d\varphi(x)} {|\nabla\varphi(x)|})$.
Итак, в результате мы имеем дуальность таких понятий как
минимальные поверхности и минимальные потоки.

\end{document}